\documentclass[a4paper,12pt]{article}
\usepackage{amssymb,latexsym,amsmath,amsthm,amscd,graphpap}
\usepackage{amsmath}
\usepackage{amsfonts}
\usepackage{amssymb}
\usepackage[dvips]{epsfig}

\newtheorem{theorem}{Theorem}[section]
\newtheorem{prop}[theorem]{Proposition}
\newtheorem{claim}[theorem]{Claim}

\newtheorem{defi}[theorem]{Definition}
\newtheorem{defis}[theorem]{Definitions}

\newtheorem{coro}[theorem]{Corollary}

\newenvironment{demo}{ \noindent \emph{\textbf{Proof:}}}{\hfill$\square$\\
\vspace{0.4cm}}
\newenvironment{sketch}{ \noindent
\emph{\textbf{Sketch of proof:}}}{\hfill$\square$\\
\vspace{0.4cm}}

\newcommand{\RR}{\mathbb{R}}

\newcommand{\TT}{\mathbb{T}}
\newcommand{\Lc}{\mathcal{L}}

\newcommand{\NN}{\mathbb{N}}

\newcommand{\Gg}{\mathfrak{G}}
\newcommand{\Xg}{\mathfrak{X}}

\newcommand{\Cc}{\mathcal{C}}
\newcommand{\Nc}{\mathcal{N}}
\newcommand{\Mc}{\mathcal{M}}

\newcommand{\Oc}{\mathcal{O}}

\newcommand{\Uc}{\mathcal{U}}

\newcommand{\ZZ}{\mathbb{Z}}

\newcommand{\no}{n$^{\text{o}}$}

\numberwithin{equation}{section}

\newdimen\texpscorrection
\newdimen\figcenter
\def\figurewithtex #1 #2 #3 #4 #5\cr{\null
  {\goodbreak\figcenter=\hsize\relax
  \advance\figcenter by -#4truecm
  \divide\figcenter by 2
  \begin{figure}[hbt]
  \vskip #3truecm\noindent\hskip\figcenter
  \includegraphics{#1}{\hskip\texpscorrection\input #2 }
  \vskip 0.8truecm{\baselineskip=0.8\baselineskip
  \noindent \vbox{\noindent {\footnotesize #5}}\par}
  \end{figure}}}
\def\point#1 #2 #3 {\rlap{\kern #1 truecm
\raise #2 truecm \hbox{#3}}}

\begin{document}

\title{\bf Observation and inverse
problems in coupled cell networks}

\author{Romain Joly}

\date{March 2011}

\maketitle
\vspace{1cm}

\begin{abstract}
A coupled cell network is a model for many situations such as food webs in
ecosystems, cellular metabolism, economical networks... It consists in a
directed graph $G$, each node (or cell) representing an agent of the network
and each directed arrow representing which agent acts on which one. It yields a
system of differential equations 
$\dot x(t)=f(x(t))$, where the component $i$ of $f$
depends only on the cells $x_j(t)$ for which the arrow $j\rightarrow i$ exists
in $G$. In this paper, we investigate the observation problems in
coupled cell networks: can one deduce the behaviour of the whole network
(oscillations, stabilisation etc.) by observing only one of the cells? We show
that the natural observation properties holds for almost all the interactions
$f$.\\[3mm]
\end{abstract}

\vspace{1cm}

\section{Introduction}

{\noindent \bf The coupled cell networks.}\\
In the recent years, the mathematical study of coupled cell networks has
been quickly developing. It combines several interests: it is strongly
related with
applications and real phenomena, the setting is very simple and it
leads to a rich class of mathematically interesting problems.
A coupled cell network models a group of agents, each one interacting with a
given part of the others through differential equations. It is represented by an
directed graph, each node being one of the agents and each directed arrow
representing which agent acts on which one. This modelling appears in many
concrete situations: networks of neurons, cellular metabolic networks, 
food webs in ecosystems, economic networks etc. Many good arguments for
studying coupled cell networks are given in \cite{Stewart}. Many examples of
coupled cell networks can also be found in \cite{GS} and \cite{SGP} and the
references therein.

The mathematical setting is the following. Let $G$ be a directed graph
with $N$ cells linked with arrows. To the cell $i$, we associate a
phase space $X_i$ of finite dimension $d_i$, which is assumed in this paper to
be the torus $(\RR / \ZZ)^{d_i}=\TT^{d_i}$ (this assumption is made for
sake of simplicity in our proofs, but our results also hold if $X_i$ is a more
general manifold, including $X_i=\RR^{d_i}$, see Section \ref{sect-further}). We
set $X=X_1\times\ldots\times X_N$ and
$d=d_1+\ldots+d_N$. For any $x\in X$ and any set of indices
$I=\{i_1,\ldots,i_k\}$, $x_I$ denotes the image of the canonical projection of
$x$ onto $X_I=X_{i_1}\times\ldots\times X_{i_k}$ (if $I=\{i\}$, we simply write
$x_i$). We also denote by $d_I=d_{i_1}+\ldots+d_{i_k}$ the dimension of $X_I$. 
For each cell $i$, the {\it direct inputs} of $i$ is the set of cells $j$ such
that the arrow $j\rightarrow i$ belongs to $G$. A node $j$ is an {\it indirect
input} of $i$ if there exists a sequence of arrows $j\rightarrow k_1 \rightarrow
k_2 \rightarrow \ldots \rightarrow i$ belonging to $G$.

Let $\Xg^1(X)$ be the space of the $\Cc^1-$vector fields on $X$ endowed with
the usual Banach topology (see for example \cite{Palis-de-Melo} or
\cite{Abraham-Robbin}). To simplify the notations, for any $x\in X$, we identify
$T_xX$ with $\RR^d$ and $\Xg^1(X)$ with $\Cc^1(X,\RR^d)$. We introduce the
subspace of all the {\it admissible vector fields}
\begin{align*}
\Cc_G&=\{f\in\Xg^1(X)~,~f_i\text{ depends
only on the direct inputs of the cell }i\}
\end{align*}
With each $f\in\Cc_G$, we associate the dynamical
system $S(t)x^0\equiv x(t)$, called {\it coupled cell network}, generated by the
differential equation
\begin{equation}\label{eq}
\left\{\begin{array}{l}
\dot x(t)=f(x(t))~,~t\in\RR\\
x(0)=x^0\in X
\end{array}\right.
\end{equation}
Notice that, since $X$ is compact and $f$ is of class $\Cc^1$, the solutions of
\eqref{eq} exist for all times.

\begin{figure}[tp]
\begin{center}
\begin{minipage}{5cm}
 \epsfig{width=4cm,file=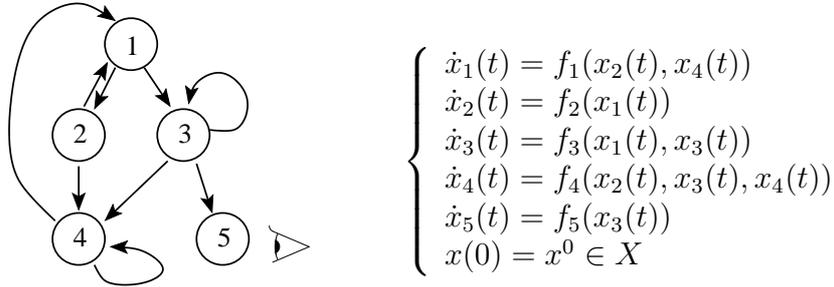}
\end{minipage}
\begin{minipage}{6cm}
$$
\left\{\begin{array}{l}
\dot x_1(t)=f_1(x_2(t),x_4(t))\\
\dot x_2(t)=f_2(x_1(t))\\
\dot x_3(t)=f_3(x_1(t),x_3(t))\\
\dot x_4(t)=f_4(x_2(t),x_3(t),x_4(t))\\
\dot x_5(t)=f_5(x_3(t))\\
x(0)=x^0\in X
\end{array}\right.
$$
\end{minipage}
\caption{\it An example of directed graph $G$ and the associated ordinary
differential equation generating the coupled cell network dynamical system. The
cell $5$ is the only observation cell of the graph $G$. The set of cells
$\{1,2,3,4\}$ is an independent strongly connected sub-network.}
\label{fig1}
\end{center}
\end{figure}

The properties of $S(t)$ may strongly depend on the structure of the directed
graph $G$. We introduce here some definitions.
\begin{defis}\label{defi-graph}
$~$

We say that a cell $i$
is an \emph{observation cell} if any cell $j$ of $G$ with $i\neq j$ is an
indirect input of $i$. 

We say that a set of cells $I$
is an \emph{independent sub-network} if $I$ contains all its indirect inputs,
i.e. if the system \eqref{eq} restricted to the set of cells $I$ is still an
autonomous ODE. 

We say that the graph $G$ is \emph{strongly connected} if any cell is an
observation cell or equivalently if any cell is an indirect input of any other
one.

We say that the graph $G$ is \emph{self-dependent} if any cell is a direct
input of itself.

We say that the graph $G$ is \emph{dimensionally decreasing} (resp.
\emph{dimensionally non-in\-crea\-sing}) if the following property holds. Let
$I$
be any non-trivial set of cells i.e. $I\neq \emptyset$ and $I\neq
\{1,\ldots,N\}$. Let $J$ be its set of direct inputs. Then, the dimension $d_J$
is larger than $d_I$ (resp. larger or equal to $d_I$).
\end{defis}

\noindent We would like to emphasise several important remarks concerning the
above setting:\\
1) We recall that we have assumed that the phase space of each cell
is a torus, when it is usually $\RR^{d_i}$ or a more general manifold,
because the compactness of $X_i$ and the triviality of its tangent
bundle simplify the statement and the proof of several of our results. In
this way, we
avoid technicalities, which are not the subject of this paper. However,
we underline that the results of this article are also valid for
$X_i=\RR^{d_i}$ or other manifolds, see Section \ref{sect-further}.\\
2) Another difference between our notations and the usual ones is that
we do not introduce a relation of equivalence between cells and arrows. This
kind of relation is introduced to model symmetries as in \cite{SGP}. For
example, it is natural to assume that several neurons of the same type interact
in the same way, that is that the functions $f_i$ corresponding to these
interactions are equal. Our purpose being to prove generic results with respect
to these $f_i$, adding constraints of symmetry increases the difficulty of the
problem. Since this paper is a first step in this subject, we have chosen to
avoid this
difficulty. However, adapting the results of this article to the presence
of symmetries could be an interesting subject of future research.\\
3) A strongly connected graph is sometimes called a transitive graph or
a path-connected graph. But, in graph theory, a transitive graph is a graph 
such that if the arrows $i\rightarrow j\rightarrow k$ exist then $i\rightarrow
k$ also exists. Thus, we decided to use the term ``strongly connected'', which
seems to be the usual one in graph theory.

\vspace{3mm}

{\noindent \bf A series of observation problems.}\\
At first sight, a coupled cell network is a simple ODE and one may wonder why
it brings new topics of research. The important characteristics of the coupled
cell networks is the splitting in several cells and the constraints due to the
structure of the interactions described by the graph $G$. This yields several 
classes of problems. For example, one can study the creation and the stability
of patterns according to the symmetries of the network. This
seems to be one of the
main classes of questions, which have been studied until now (see \cite{GNS},
\cite{GS} and \cite{SGP} as well as the many related works). Another class is
the observability problems, which are the subject of this paper and has been
previously tackled in \cite{GRW}. Basically, one wonders if one can deduce
the behaviour of the whole network only by observing one of the cells. We
investigate here three problems.

\noindent {\it Problem 1: inverse problem.}\\ 
Can one distinguish two trajectories by only looking at one cell?\\
\noindent {\it Problem 2: observation of oscillations.}\\ 
One of the cells is
oscillating with a period $T$, are the other cells also oscillating?\\
\noindent {\it Problem 3: observation of stabilisation.}\\ 
One of the cells is stabilising and converging to an equilibrium state, have
the other cells the same behaviour?
\begin{figure}[tp]
\begin{center}
\begin{minipage}{7cm}
 \epsfig{width=5cm,file=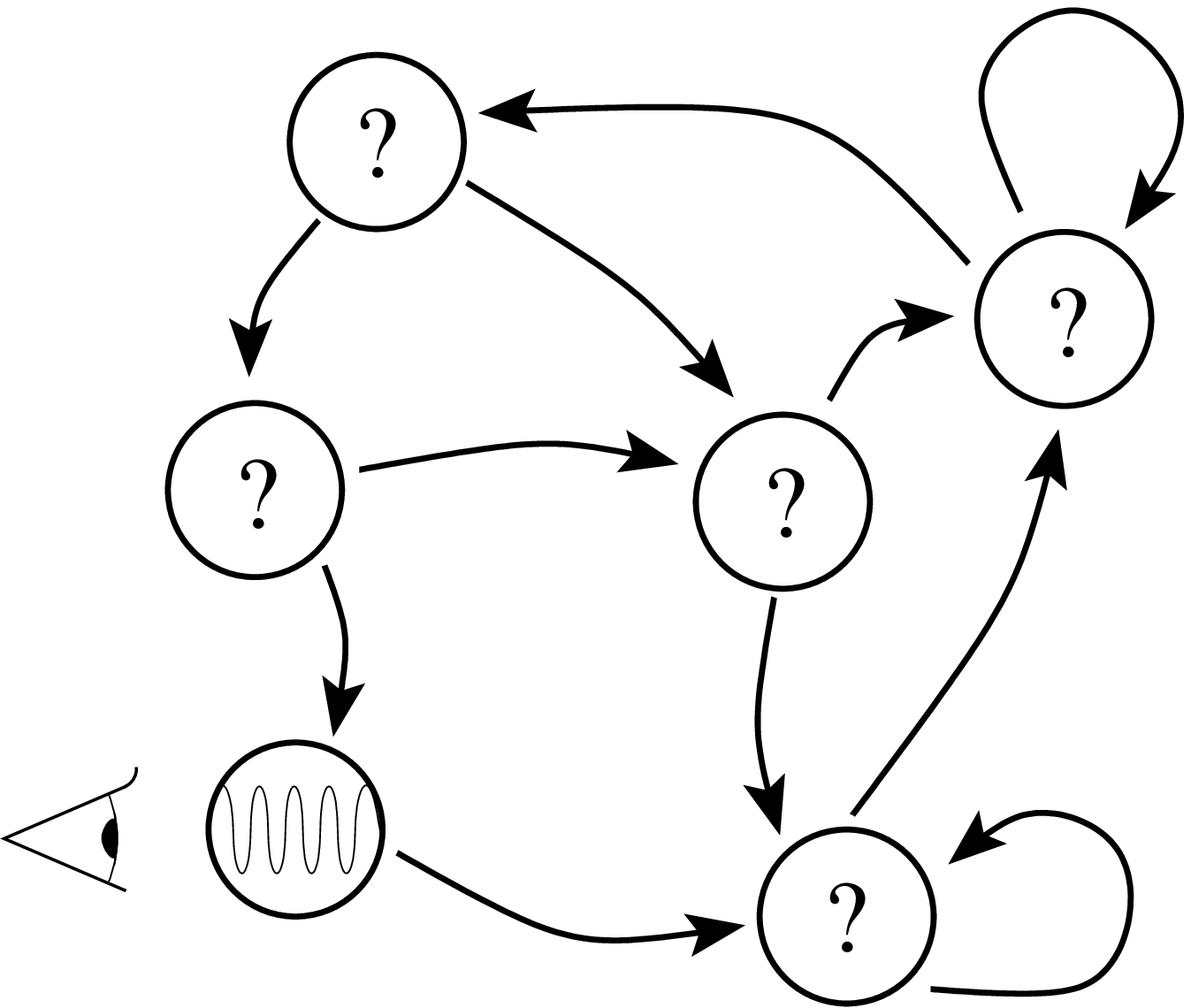}
\end{minipage}
\begin{minipage}{7cm}
 \epsfig{width=5cm,file=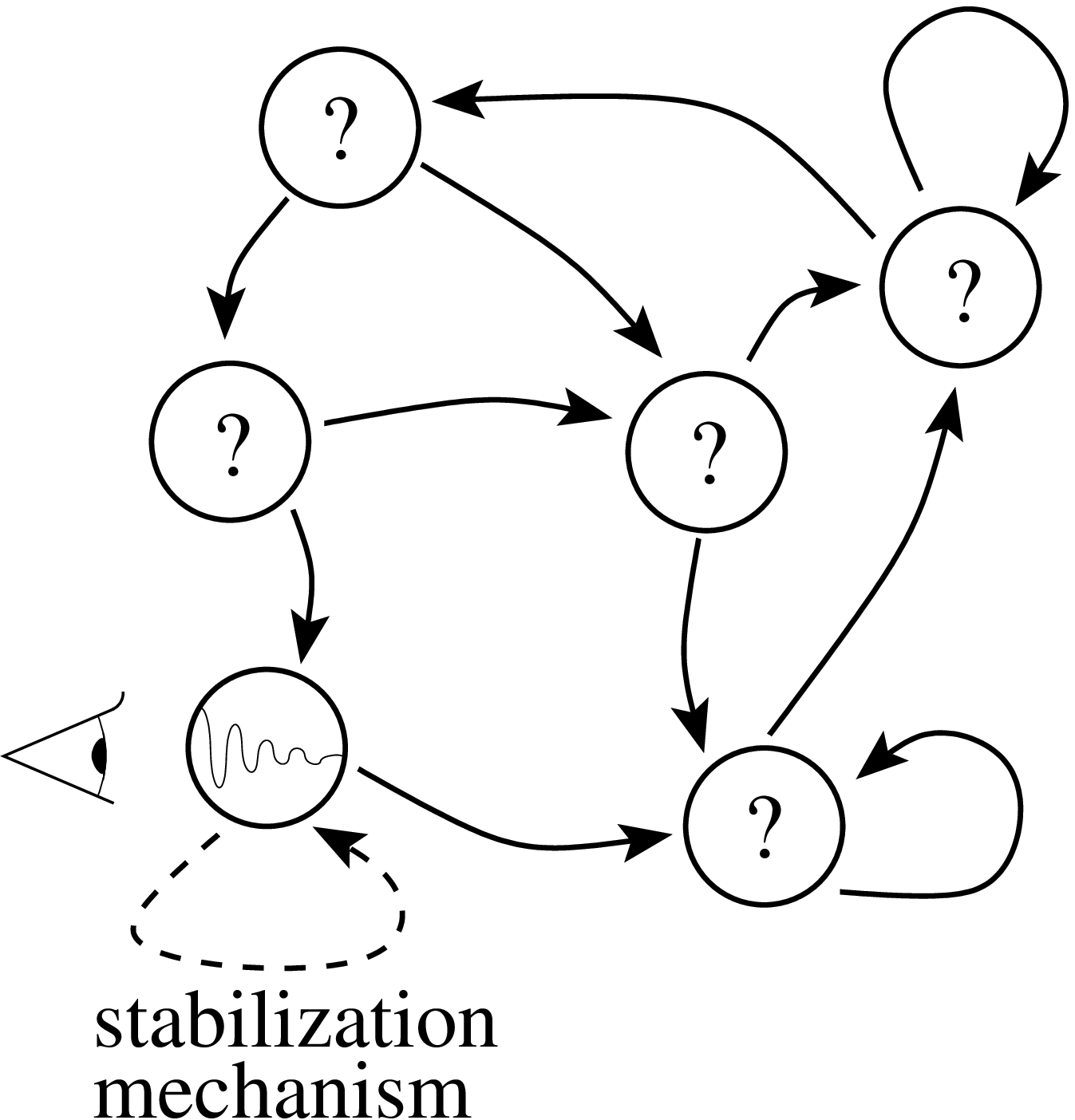}
\end{minipage}
\caption{\it Two examples of observation problems. Left, one measures
the state of one cell and observes oscillation. Do the states of the other
cells also oscillate? Right, assume that one has found a feedback mechanism
that stabilises the observation cell. Do the other cells also stop oscillating?}
\label{fig2}
\end{center}
\end{figure}

These kinds of questions are natural in the context of coupled cell
network modelling. For example, if one considers an economical network, it is
impossible to measure the behaviour of all the economical agents. One would
like to know if the measures of a small number of indices is sufficient to know
if the economy is chaotic, oscillating or stable. Inverse problem are also
motivated by the desire of controlling a whole network that can be
observed only through a small number of parameters (as it is the case for a
whole organism or an ecosystem). If one cannot
distinguish two trajectories by observing a small part of the network, then one
cannot know in which state the whole network is and therefore what control one
should impose.

\vspace{3mm}

{\noindent \bf The necessity and the meaning of generic results.}\\
The first answer to problems 1-3 is simply ``no''. Of course, there is firstly 
a problem with the geometry of the graph $G$: two cells have totally
independent behaviours if they are not interacting, even indirectly. But, even
with a nice strongly connected network, it is not difficult
to construct particular interactions $f_i$ such that the observation properties
fail.

So, at a second time, one could be tempted to obtain a criterion on $f$, for
which observability is possible. However, what is the usefulness of a criterion
such as ``these two biological parameters must have irrational ratio''? In
physics,
interactions have a well-known structure. But, in biology or economy, one may
have no precise idea of the form of the interactions $f_i$. In these cases, a
concrete criterion has no meaning.

In this situation, one must use the concept of generic properties. The idea is
to show that for almost all the possible interactions, problems 1-3 have
positive answers. Then, one believes that, when one deals with a concrete
network,
the observations properties surely hold for this network. In this paper, we
fix a given directed graph $G$ and by ``almost all the possible interactions'',
we mean a generic subset of the admissible vector fields $\Cc_G$. We recall
that a generic subset is a subset that contains a countable intersection of
dense open subsets. Since $\Cc_G$ is a Banach space, generic subsets are dense
and the use of the notion of generic subsets to characterise large subsets is
well admitted. 

Generic results are meaningful in network's problems, where the functions of
interactions $f_i$ are not completely a priori determined (typically in food
webs in ecosystems, in metabolic networks...), so that one can hope
that the considered network is a generic one. If $f$ is precisely determined
(typically in networks following simple physical laws), a generic result is
meaningless. If the interactions follow some symmetries (typically in neuronal
networks), then one should prove a generic result in the class $\Cc_G$
restricted to interactions having these symmetries.

\vspace{3mm}

{\noindent \bf Main results}\\
The inverse problem, Problem 1, depends on the geometry of the network
$G$ (it is noteworthy that it is our only theorem which does not hold for any
strongly connected graph). The inverse problem may fail if there are some cells
on which the trajectories are constant (see the counter-example of Section
\ref{sect-eq-diff}). In Sections \ref{sect-eq} and \ref{sect-inv}, we state 
different properties related to this inverse problem, which give a sharp
idea of the cases where it may not hold. They may be summarised in the following
simple statement.
\begin{theorem}\label{th1}
{\bf Observation of trajectories}\\
Let $G$ be a given directed graph.
There exists a generic subset $\Gg^{traj}$ of $\Cc_G$ such that for any
$f\in\Gg^{traj}$ the following property holds. Let $x(t)$ and $y(t)$ be two
solutions of \eqref{eq} and let $i$ be an observation cell. Assume that there is
a time interval $(a,b)$, $a<b$, such that $x_i(t)=y_i(t)$ for all $t\in(a,b)$.
Then, the
whole trajectories $x(t)$ and $y(t)$ coincide, except maybe in the cells where
the trajectories $x(t)$ and $y(t)$ are both constant. Moreover, if the graph $G$
is either a self-dependent directed graph, or a dimensionally decreasing graph,
then $x(t)$ and
$y(t)$ coincide everywhere, even in the cells where the trajectories $x(t)$ and
$y(t)$ are both constant.
\end{theorem}

The observation of oscillations stated in problem 2 is a natural problem,
which has already been studied in \cite{GRW}. It is proved there that, in
self-dependent and strongly connected graph, hyperbolic periodic orbits have
generically no constant cells. We obtain here a stronger result in several
directions: no hyperbolicity of periodic orbit is needed, no self-dependence is
assumed, all the cells are not only oscillating but also oscillate with the same
period and the constant cells are observable only in constant solutions.

In Theorem \ref{th2} below, we use the following notations. We say that a
function $p(t)$ is {\it $T-$periodic} on $(a,b)$ if $p(t)=p(t+T)$ for any $t\in
(a,b-T)$. We say that a function $p(t)$ is {\it exactly $T-$periodic} on
$(a,b)$ if $T> 0$ is the smallest number such that $p$ is $T-$periodic on
$(a,b)$. 

\begin{theorem}\label{th2}
{\bf Observation of oscillations}\\
Let $G$ be a given directed graph.
There exists a generic subset $\Gg^{osc}$ of $\Cc_G$ such that for any
$f\in\Gg^{osc}$ the following properties hold. Let $x(t)$ be a solution of
\eqref{eq}, then
\begin{itemize}
\item[(a)] if there is a cell $i$ such that $x_i(t)$ is constant on
some time interval $(a,b)$, then $x_i(t)$, as well as any $x_j(t)$ with
$j$ being an indirect input of the cell $i$, is constant on $\RR$. In
particular, if $i$ is an observation cell, then $x(t)$ is an equilibrium point
of \eqref{eq}.
\item[(b)] if there is a cell $i$, a positive time $\tau$ and a positive period
$T$ such that $x_i(t)$
is $T-$periodic on some time interval $(a,a+T+\tau)$, then $x_i(t)$, as well as
any $x_j(t)$ with $j$ being an indirect input of the cell $i$, are
$T-$periodic in $\RR$.
\item[(c)] if $G$ is strongly connected and if there is a cell $i$ where
$x_i(t)$ is exactly $T-$periodic on some time interval $(a,a+T+\tau)$, then the
whole orbit $x(t)$ and any of the $x_j(t)$ are exactly $T-$periodic.
\end{itemize}
\end{theorem}
Notice that Theorem \ref{th2} is not in contradiction with the examples of
\cite{GNS}, even the ones which are stable with respect to perturbations of $f$.
Indeed, we assume here no a priori symmetry for the interactions $f_i$.

Our last result concerns Problem 3, that is the observation of stabilisation.
It is a natural problem since the stabilisation of a network is an
important issue in many applications: economy, ecosystems, metabolism...
\begin{theorem}\label{th3}
{\bf Observation of stabilisation}\\
Let $G$ be a given directed graph.
There exists a generic subset $\Gg^{stab}$ of $\Cc_G$ such that for any
$f\in\Gg^{stab}$, the following property holds. Let $x(t)$ be a
solution of \eqref{eq} and let $i$ be an observation cell. If $x_i(t)$
converges to a constant $x^*_i$ when $t$ goes to $+\infty$, then the whole
trajectory $x(t)$ converges to a single equilibrium point $x^*\in X$ when $t$
goes to $+\infty$.
\end{theorem}
We underline that, due to Theorem \ref{th1}, the equilibrium $x^*$ can be
moreover uniquely determined by observing $x^*_i$ if $G$ is
either a self-dependent directed graph or a dimensionally decreasing graph.

\vspace{3mm}

{\noindent \bf Organisation of the paper.}
The main tools for proving Theorems \ref{th1}, \ref{th2} and \ref{th3} are
transversality theorems. We recall the ones that we need in Section
\ref{sect-pre}. In Section \ref{sect-eq}, we investigate the observation
problem for equilibrium states of the network. We also prove the generic
simplicity of equilibrium points, which is an interesting property by itself.
The inverse problem
stated in Theorem \ref{th1} is proved in Section \ref{sect-inv} and, in Section
\ref{sect-coro}, Theorem \ref{th2} and \ref{th3} are deduced from the previous
results. Finally, we present further results in Section \ref{sect-further}, in
particular the non-compact case of $X_i=\RR^{d_i}$ and the case of discrete
networks.

\vspace{5mm}

{\noindent \bf Acknowledgements: } The author would like to thank Martin
Golubitsky and Genevi\`eve Raugel for pointing to him the interest of coupled
cell networks. Several discussions have also been useful for the writing of
this paper: Matthieu L\'eautaud and Vincent Perrollaz have enhanced to the
author the interest of studying Problem 1 and 3 and Beno\^\i t Kloeckner has
patiently answered many questions about differential geometry.


\section{Geometric preliminaries}\label{sect-pre}
\subsection{Transversality theorems}\label{sect-trans}
In a space $X$, which admits a natural measure, the notion of ``almost
everywhere'' as sets of full measure is natural and well admitted. As soon as
the space $X$ is infinite-dimensional and has no natural measure, the
meaning of ``almost everywhere'' is less clear. If $X$ is a Baire space
(typically a complete space), then the most used notion of ``almost
everywhere'' is the notion of generic subsets, that are subsets of $X$
containing a countable intersection of dense open subsets of $X$. A property is
said to be generic in $X$ (or with respect to $x\in X$) if it holds for all $x$
in a generic subset of $X$. It is noteworthy that there are other convenient
notions of ``almost everywhere'' in Banach spaces such as the notion of
prevalence, see Section \ref{sect-further}.

In problems involving finite-dimensional manifolds, the proofs of the genericity
of a property mainly use Sard's Theorem or theorems of transversality similar to
the ones of Thom (see \cite{Sard}, \cite{Thom1} and \cite{Thom2}). We recall
here the classical Sard-Smale Theorem (see for example \cite{Abraham-Robbin} for
a proof).
\begin{defi}\label{defi-reg}
 Let $\Mc$ and $\Nc$ be two $\Cc^1-$ manifolds and let $f:\Mc\rightarrow\Nc$ be
a $\Cc^1-$map. A point $x\in\Mc$ is \emph{regular point} of $f$ if $Df(x)$ is
surjective, otherwise $x$ is a \emph{critical point} of $f$. A point
$y\in\Nc$ is a \emph{regular value} of $f$ if any $x\in f^{-1}(y)$ is a
regular point of $f$, otherwise $y$ is a \emph{critical value} of $f$.
\end{defi}
\begin{theorem}\label{th-Sard-Smale}
{\bf A transversality Theorem}\\
Let $r\geq 1$. Let $\Mc$ and $\Nc$ be two finite-dimensional $\Cc^r-$manifolds
and let $\Lambda$ be a $\Cc^r$ Banach manifold.
Let $y\in \Nc$ and let $\Phi\in \Cc^r(\Mc\times\Lambda,\Nc)$. Assume that:\\
(i)  $r>dim(\Mc)-dim(\Nc)$,\\
(ii)  for any $(x,\lambda)\in \Phi^{-1}(\{y\})$, $D\Phi(x,\lambda):T_x\Mc\times
T_\lambda \Lambda \rightarrow T_y\Nc$ is surjective,\\
(iii) $\Mc$ is separable.\\ 
Then, the set $\Gg=\{\lambda\in\Lambda,~y\text{ is a regular value of the map
}x \mapsto \Phi(x,\lambda)\}$ is a generic subset of
$\Lambda$.\\
Moreover, if $dim(\Mc)<dim(\Nc)$, then the set $\Gg$ is equal to the set
$\{\lambda\in\Lambda,~y$ is not in the image of the map $x
\mapsto \Phi(x,\lambda)\}$.
\end{theorem}

Smale in \cite{Smale} showed that Sard's Theorem can be extended to Banach
spaces by using the
notion of Fredholm operators. Later, Quinn in \cite{Quinn} noticed that
the notion of left-Fredholm operators is often sufficient. 
\begin{defi}\label{defi-Fred}
Let $X$ and $Y$ be two Banach spaces. A bounded linear operator $L:X\rightarrow
Y$ is a \emph{left-Fredholm operator} if:\\
(i) its kernel $Ker(L)$ splits in $X$, i.e. there exists a space $X_1$ such
that $X=X_1\oplus Ker(L)$ and $X_1$ and $Ker(L)$ are closed subspaces,\\
(ii) its image $R(L)$ splits in $Y$, i.e. there exists a space $Y_2$ such
that $Y=R(L)\oplus Y_2$ and both subspaces are closed,\\
(iii) its kernel $Ker(L)$ is finite-dimensional.\\
If moreover the supplementary space $Y_2$ is also finite-dimensional, then $L$
is called a \emph {Fredholm operator}.\\
The index of $L$ is defined by $Ind(L)=dim(Ker(L))-dim(R(L))$ (which is equal to
$-\infty$ if $L$ is not a Fredholm operator).
\end{defi}
Following Smale arguments, one can extend the classical transversality theorems,
as Theorem \ref{th-Sard-Smale}, to Banach manifolds. There exist many different
versions of this kind of theorems in Banach manifolds (often called Sard-Smale
theorems), see for example \cite{Abraham-Robbin},
\cite{Henry} or \cite{Saut-Temam}. We will use here a version proved by Henry in
\cite{Henry} (see also \cite{RJ-prevalence}, where its extension to the
notion of prevalence is proved).
\begin{theorem}\label{th-Henry}
{\bf Henry's Theorem}\\
Let $\Mc$, $\Lambda$ and $\Nc$ be three Banach manifolds. Let $\Phi:\Mc\times
\Lambda \longrightarrow \Nc$ be a map of class $\Cc^1$ and $y$ be a point of
$\Nc$.\\ 
We assume that :
\begin{description}
\item{(i)} $\forall (x,\lambda)\in \Phi^{-1}(y)$, 
$D_x\Phi(x,\lambda):T_x\Mc \rightarrow T_y(\Nc)$ is a left-Fredholm 
  operator with negative index,
\item{(ii)} $\forall (x,\lambda)\in \Phi^{-1}(y)$,  the image of the total
derivative $D\Phi(x,\lambda):T_x\Nc\times
T_\lambda\Lambda\rightarrow T_y\Nc$ contains a finite-dimensional subspace $Z$
such that $Z\cap R(D_x\Phi(x,\lambda))=\{0\}$ and the dimension of $Z$ is
strictly larger than the one of $Ker(D_x\Phi(x,\lambda))$,
\item{(iii)} $\Mc\times\Lambda$ is separable.
\end{description}
Then there exists a generic subset $\Gg$ of $\Lambda$ such
that, for any $\lambda_0 \in\Gg$, $y$ is not in the image of the map $x\mapsto
\Phi(x,\lambda_0)$. If moreover $\Mc$ is compact, then $\Gg$ is an open set.
\end{theorem}

In our applications of Theorem \ref{th-Henry}, the operator $D_x\Phi$ can be
split as $D_x\Phi=L+K$, where $K$ a compact operator and $L$ is a simple
operator, for which Hypotheses (i) and (ii) are easily checked. Therefore,
the following propositions will be useful.
\begin{prop}\label{prop-Fred-1}
Let $X$ and $Y$ be two Banach spaces. Let $L:X\rightarrow Y$ be a
left-Fredholm  operator and let $K\in\Lc(X,Y)$ be a compact operator. Then,
$L+K$ is a  left-Fredholm map with the same index as the one
of $L$. 
\end{prop}
The proof of Proposition \ref{prop-Fred-1} is classical, see for
example \cite{Bonic}.
To check Hypothesis (ii) of Theorem \ref{th-Henry}, we will use the following
criterion.
\begin{prop}\label{prop-Fred-3}
Let $L:X\rightarrow Y$ be a left-Fredholm map and $X_1$ and $Y_2$ be as in
Definition \ref{defi-Fred}. Assume that $(z_n)_{n\in\NN}$ is 
an infinite free family of $Y_2$. Then, for any compact operator $K:X\rightarrow
Y$, there exists a finite dimensional subspace
$Z=span(z_{n_1},\ldots,z_{n_p})$ such that $Z\cap R(L+K)=\{0\}$
and the dimension of $Z$ is strictly larger than the one of
$Ker(L+K)$.
\end{prop}
\begin{demo}
 Let $P$ be the continuous projection on $R(L)$ canonically defined by the
closed splitting $Y=R(L)+Y_2$. Let $K$ be a given compact operator. Then, $PK$ 
is a compact operator and thus Proposition \ref{prop-Fred-1}
shows that $L+PK$ is a left-Fredholm operator. In particular, $L+PK$ has a
finite-dimensional kernel. On the other hand, $(L+PK)x=0$ if and only if 
$P(L+K)x=0$ that is that $(L+K)x$ belongs to $Y_2$. Thus, $dim(R(L+K)\cap
Y_2)=dim(Ker(L+PK))<\infty$ and an at most finite
dimensional subspace of $span(z_n)$ is included in $R(L+K)$. Hence, we can
still extract a subspace $Z$ from $span(z_n)$, with dimension as large as
needed, and such that $Z\cap R(L+K)=\{0\}$.
\end{demo}

Let us notice that the assumptions that the sequence $(z_n)$ belongs to a closed
complementary space of $R(L)$ and that $L$ is left-Fredholm are important, as
shown by the following examples. We set $X=Y=\ell^\infty(\RR)$ and we define
$e_i\in Y$ by $e_i(n)=\delta_{i=n}$, i.e. $(e_i)_{i\in\NN}$ is the canonical
``basis''. Proposition \ref{prop-Fred-3} does not hold for
$$ \left\{\begin{array}{l} L(x_1,x_2,x_3,\ldots)=
(x_1,0,x_2,0,x_3,\ldots)\\
K(x_1,x_2,x_3,\ldots)=(0,x_1,0,\frac{x_2}2,0,\frac{x_3}3,\ldots)\\
z_n=e_{2n-1}+\frac {e_{2n}}{n} \end{array}\right.
\text{ or }
\left\{\begin{array}{l} 
L(x_1,x_2,x_3,\ldots)=(x_1,0,x_3,0,x_5,\ldots)\\
K(x_1,x_2,x_3,\ldots)=(0,\frac{x_2}2,0,\frac{x_4}4,0,\frac{x_6}6,\ldots)\\
z_n=e_{2n}\end{array}\right.
$$
In the first example, $span(z_n)$ is not in a closed complementary subset of
$R(L)$. In the second one, $L$ is not left-Fredholm.

\subsection{The Banach manifold of $\Cc^1-$paths}
In this section, we would like to discuss the different topologies of the spaces
used in this paper. First, we embed $\RR^d$ with the supremum norm
$\|.\|_\infty$ and we identify the torus $X=\TT^d$ to $(\RR/\ZZ)^d$
with its quotient topology. In
particular, the distance $d(x,y)$ between two points $x$ and $y$ of $\TT^d$ is
bounded by $1/2$. We also recall that we have identified all the tangent
spaces of $X$ to $\RR^d$.

In this paper, the space of parameters $\Lambda$ of Theorems
\ref{th-Sard-Smale} and \ref{th-Henry} will always be the space of admissible
vector fields $\Cc_G$, which is identified to a closed
subspace of $\Cc^1(X,\RR^d)$ with its natural Banach structure.

In some of our proofs, the Banach manifold $\Mc$ of
Theorem \ref{th-Henry} will be the space $\Cc^1([0,\tau],X)$ of $\Cc^1-$paths on
$X$, with $\tau>0$ a given time. Its structure of Banach manifold is as follows.
Let $\gamma$ be a given $\Cc^1-$path. We consider a neighbourhood of $\gamma$
given by
$$\Oc=\{\beta\in\Cc^1([0,\tau],X)~,~\forall~t\in[0,\tau]~,
~d(\beta(t),\gamma(t))<1/4~\}~.$$
We now use a canonical lifting $\varphi$ along the path $\gamma$, mapping $\Oc$
to an open subset $\Uc$ of $\Cc^1(\RR^d,\RR^d)$. Notice that, by definition of
$\Oc$, lifting the whole neighbourhood $\Oc$ is possible and moreover
\begin{equation}\label{c1path}
|\varphi(\gamma_i(t))-\varphi(\beta_i(t))|=d(\gamma_i(t),\beta_i(t))
\end{equation}
for any $\beta\in\Oc$, any $t\in[0,\tau]$ and any cell $i$. Therefore, we can
identify $\Oc$ to $\Uc$, which is an open subset of a Banach space. Notice
that the tangent space to $\Cc^1([0,\tau],X)$ at a path $\gamma$ is simply
$\Cc^1([0,\tau],\RR^d)$.


\section{Generic properties of equilibrium points}\label{sect-eq}

In this section, we prove several generic properties of equilibria of the ODE,
i.e. singularities of the vector field $f$, that are points $x\in X$ such that
$f(x)=0$. These properties are interesting by themselves, but they are also
important steps in proving Theorem \ref{th1} and \ref{th2}.

\subsection{Generic simplicity of equilibria}\label{sect-simpl}
We say that an equilibrium point $x\in X$ is simple if the differential $Df(x)$
is surjective, which is equivalent to the injectivity of this differential. A
simple equilibrium point $x\in X$ is isolated from the other equilibria, and
hence there is at most a finite number of simple equilibria. Moreover, by the
implicit functions theorem, simple equilibria, which are non-degenerate
singularities of the vector fields, depend smoothly of small perturbations of
the vector field. That is why, the simplicity of equilibria of an ODE is an
important property to consider.
\begin{prop}\label{prop-simpl}
There exists an open dense set $\Gg^{simpl}$ of $\Cc_G$ such
that for any $f\in\Gg^{simpl}$, all the equilibrium points of \eqref{eq} are
simple.
\end{prop}
\begin{demo}
The arguments of this proof are very classical and do not differ from the ones
used for general ODEs, i.e. without the constraint that $f\in\Cc_G$. 
First, let us notice that $f\in \Gg^{simpl}$ is equivalent to $0$
being a regular value of $f$. 

The set $\Gg^{simpl}$ is open. Indeed, let $(f_n)$
be a sequence of vector fields converging to $f$ and such that there exists
$x_n\in X$ with $f_n(x_n)=0$ and $Df_n(x_n)$ is not surjective (i.e.
$f_n\not\in\Gg^{simpl}$). Since $X$ is compact, up to extracting a subsequence,
one can assume that $(x_n)$ converges to $x\in X$. By continuity, we get that
$x$ is both an equilibrium point and a critical point of $f$, which shows that
$f\not\in\Gg^{simpl}$.

Let us show that $\Gg^{simpl}$ is dense. We set $\Mc=X$, $\Nc=\RR^d$
and $\Lambda=\Cc_G$. We consider the point $y=0\in\RR^d$ and the function
$$\Phi:~(x,f)\in\Mc\times\Lambda~\longmapsto~f(x)\in\RR^d~.$$
The function $\Phi$ is of class $\Cc^1$ and 
$$ D\Phi(x,f).(\xi,g)=D_x\Phi(x,f).\xi+D_f\Phi(x,f).g=Df(x).\xi+g(x)~.$$
Since $D_x\Phi(x,f)=Df(x)$, $y=0$ is a regular value of
$x\mapsto\Phi(x,f)$ if and only if $0$ is a
regular value of $f$, i.e. $f\in\Gg^{simpl}$. Thus, it is
sufficient to check the hypotheses of Theorem \ref{th-Sard-Smale} to get the
genericity of $\Gg^{simpl}$. We have $dim(\Mc)-dim(\Nc)=0$, thus Hypothesis
(i) holds. Hypothesis (iii) is satisfied since $X$ is compact. To check
Hypothesis (ii), it is sufficient to notice that, for any $x\in X$, the map
$g\in\Cc_G \mapsto D_f\Phi(x,f).g=g(x)$ is surjective onto $\RR^d$ since any
constant vector field $g\in\RR^d$ is an admissible vector field.
\end{demo}

\noindent {\bf Nota Bene:} the hyperbolicity of the equilibrium points is a
stronger property than simplicity. Hyperbolicity is an important concept for
studying the behaviour of dynamical systems. Even, if we do not need it to
prove our main results, we discuss the generic hyperbolicity of the
equilibrium points of \eqref{eq} in Section \ref{sect-further}.

\vspace{5mm}

For some geometries of the network $G$, the differential equation \eqref{eq}
cannot admit simple equilibrium points. Thus, in these cases, Proposition
\ref{prop-simpl} implies that \eqref{eq} has in fact no equilibrium point.
\begin{coro}\label{coro-eq}
Assume that there exists a set $I$ of cells, with a set $J$ of direct inputs,
such that $d_J<d_I$ (this holds in particular if the graph $G$ is not
dimensionally non-increasing). Then, for any $f$ belonging to $\Gg^{simpl}$, the
open subset of $\Cc_G$ introduced in Proposition \ref{prop-simpl}, the
differential equation \eqref{eq} has no equilibrium point. 
\end{coro}
\begin{demo}
Assume that $G$ is not dimensionally non-increasing and let $I$ and $J$ be such
that $d_J<d_I$ and $J$ is the set of the direct inputs of the cells of 
$I$. If $x$ is a simple equilibrium point of \eqref{eq}, the
differential $Df(x)$ should be surjective, thus this should be a fortiori the
case for the restricted map $\xi\in \RR^d \mapsto Df_I(x).\xi \in \RR^{d_I}$, or
more precisely, omitting the unnecessary cells, the map $\xi_J\in \RR^{d_J}
\mapsto Df_I(x).\xi_J\in \RR^{d_I}$. Obviously, this could not happen since
$d_J<d_I$. Therefore, \eqref{eq} 
cannot admit simple equilibrium points and thus if $f\in\Gg^{simpl}$, then
\eqref{eq} does not have any equilibrium.
\end{demo}

Corollary \ref{coro-eq} applies in particular when a cell $i$ has no input
(take $I=\{i\}$) or when a cell $i$ is the input of no other cell (take
$I$ being the whole set of cells).

\subsection{Generic observability of equilibrium points}\label{sect-eq-obs}

In this section, we prove property (a) of Theorem \ref{th2}, which is a direct
consequence of the following property. 

\begin{prop}\label{prop-eq}
Let $i$ and $j$ be two cells such that $j$ is a direct input of $i$. Let
$\tau>0$ be given and let $\Gg^{eq}_{i,j,\tau}$ be the set of all the
admissible vector fields $f\in\Cc_G$ such that, if $x(t)$ is a solution of
\eqref{eq} such that $x_i(t)$ is
constant on the time interval $[0,\tau]$, then $x_j(t)$ is also constant on
$[0,\tau]$. Then, $\Gg^{eq}_{i,j,\tau}$ is a generic subset of $\Cc_G$.
\end{prop}
\begin{demo}
Let $\tau>0$, let $i$ be a cell and let $j\neq i$ be a direct input of $i$.
We apply Theorem
\ref{th-Henry} to the following setting. The set of parameters is
$\Lambda=\Cc_G$, the space of the admissible vector fields. We set  
$\Nc=\Cc^0([0,\tau],\RR^d)$, $y=0$ and 
$$\Mc=\{\gamma\in\Cc^1([0,\tau],X),~\gamma_i(t)\text{ is constant and
}\gamma_j(t)\text{ is not constant }\}~.$$
We define the function $\Phi\in\Cc^1(\Mc\times\Lambda,\Nc)$ by
$$\Phi(\gamma,f)=\frac{d}{dt}\gamma(t)-f(\gamma(t))~.$$
The tangent space $T_\gamma\Mc$ is given by the functions
$\omega\in\Cc^1([0,\tau],\RR^d)$ such that $\omega_i(t)$ is constant.
We have
$$D\Phi(\gamma,f).(\omega,g)=\frac{d}{dt}
\omega(t)-Df(\gamma(t)).\omega(t)- g(\gamma(t))~.$$
{\noindent\it Step 1: $D_\gamma\Phi$ is a left-Fredholm map.}\\
The map $\omega \in T_\gamma\Mc \mapsto \frac{d}{dt}
\omega \in \Cc^0([0,\tau],\RR^d)$ is a left-Fredholm function. Indeed, its
kernel is the set of the constant functions and therefore is of dimension $d$
and admits a closed complementary set consisting in the functions of null
integral.
Moreover, its image is the set of functions $y(t)\in\Cc^0([0,\tau],\RR^d)$ such
that $y_i(t)\equiv 0$, and therefore is closed, and it admits
an infinite-dimensional closed complementary space $Y_2$, which is the
space of
functions $y\in\Cc^0([0,\tau],\RR^d)$ such that $y_k(t)=0$ if $k\neq i$. Using
Ascoli's theorem, we show that the injection
$\Cc^1([0,\tau],\RR^d)$ in $\Cc^0([0,\tau],\RR^d)$ is compact. 
Since $Df(\gamma(t))$ belongs to $\Cc^0([0,\tau],\Lc(\RR^d))$,
we get that the map $\omega \in
\Cc^1([0,\tau],\RR^d)\mapsto Df(\gamma(t)).\omega\in\Cc^0([0,\tau],\RR^d)$ is
compact. Therefore, using Proposition \ref{prop-Fred-1}, we obtain that
$D_\gamma\Phi$ is a left-Fredholm map.

{\noindent\it Step 2: checking Hypothesis (ii) of Theorem \ref{th-Henry}}\\
We have seen that $D_\gamma\Phi$ can be written $L+K$ as in Proposition
\ref{prop-Fred-3}, $Y_2$ being 
the space of functions $y\in\Cc^0([0,\tau],\RR^d)$ such that $y_k(t)=0$ if
$k\neq
i$. Using Proposition \ref{prop-Fred-3}, we notice that it is sufficient to 
exhibit a sequence of functions $(g^n)\in \Cc_G$ such that
$D_f\Phi.g^n=-g^n(\gamma(t))$ compose a free family of $Y_2$.
By assumption, the curve $\gamma_j$ is not reduced
to a point of $X_j$. One can choose an infinite number of disjoint open sets
$\Oc_j^n\subset X_j$ such that the curve $t\mapsto \gamma_j(t)$
belongs to $\Oc_j^n$ at some time $t^n\in[0,\tau]$. Let
$\chi^n\in\Cc^1(X_j,\RR)$ be bump functions with compact support in
$\Oc_j^n$ and let $\zeta\in\RR^d$ be such that $\zeta_k=0$ if $k\neq i$ and
$\zeta_i\neq 0$. Then, the vector fields $g^n(x)=\zeta\chi^n(x_j)$ belong to
$\Cc_G$ and are
such that the functions $t\mapsto g^n(\gamma(t))$ belong to $Y_2$ and have
disjoints supports. This yields the required free family of $Y_2$ and shows
that Hypothesis (ii) holds.

{\noindent\it Step 3: $\Mc$ and $\Lambda$ are obviously separable}

{\noindent\it Step 4: conclusion}\\
We can apply Theorem \ref{th-Henry}, which shows that there exists a generic
subset of admissible vector fields $f\in\Cc_G$ such that $0$ is not in the
image of $\Phi(.,f)$. This means that there is no solution
$\gamma\in\Cc^1([0,\tau],X)$ of
$$\frac{d}{dt}\gamma(t)-f(\gamma(t))=0$$
such that $\gamma_i$ is constant and $\gamma_j$ not constant, which is exactly
the statement of Proposition \ref{prop-eq}.
\end{demo}

{\noindent\bf Proof of Property (a) of Theorem \ref{th2}: }
We set $\Gg^{eq}=\cap_{i,j,n} \Gg^{eq}_{i,j,1/n}$,
which is a generic subset of $\Cc_G$, as a countable intersection of generic
subsets. Then, assume that $f\in\Gg^{eq}$. If
$i$ is such that $x_i(t)$ is constant on some interval $(a,b)$, then,
up to translating the time and up to reducing the time interval, one can assume
that there exists $n$ such that $x_i(t)$ is constant for $t\in
[0,1/n]$. By applying the first part of Proposition \ref{prop-eq}, $x_j(t)$ must
be constant on  $[0,1/n]$ for any $j$ direct input of $i$. Then, by recursion,
we
get that, on the interval $[0,1/n]$, any indirect input of $i$ is
constant. The set $I$ of cells containing $i$ and all its
indirect inputs is an independent sub-network and thus follows an evolution
which depends only on the values of the cells of
$I$. Therefore, the differential equation \eqref{eq} restricted to
$X_{I}$ is an ODE for which the solution $x_{I}(t)$ is constant
on $[0,1/n]$ and hence is constant for all times $t\in\RR$.
\hfill $\square$

\subsection{Generic inverse problem for equilibrium points}\label{sect-eq-diff}

One may hope that, at least generically and for strongly connected graphs, if
two equilibria coincide in one cell, then they must be equal everywhere. In
fact, a surprisingly simple counter-example exists. We look to the
circular two-cell network
\begin{equation}\label{ce-eq}
\dot x_1(t)=f_1(x_2(t))~~~~\dot x_2(t)=f_2(x_1(t))~.
\end{equation}
Equilibria consist in couples $(x_1,x_2)$ with $x_1$ and $x_2$
being singularities of the
respective vector fields $f_2$ and $f_1$. Therefore, as soon as $f_2$ has two
singularities $x_1\neq \tilde x_1$ and $f_1$ has one singularity $x_2$, then
$(x_1,x_2)$ and $(\tilde x_1,x_2)$ are both equilibria, which cannot be
distinguished by observing
the cell $2$. In addition, if the phase spaces of both cells have the same
dimension, then these singularities may be simple and thus stable with respect
to small perturbations of $f_1$ and $f_2$. Therefore, this type of
counter-example exists in an open set of non-linearities.

However, the generic inverse problem holds in some situations, as shown in the
next proposition.
\begin{prop}\label{prop-eq-diff}
Let $G$ be either\\
(a) a self-dependent directed graph,\\
(b) or a dimensionally decreasing graph,\\
(c) or a not dimensionally non-increasing graph.\\
Let $i_0$ be an observation cell for the graph $G$. Then, there exists a generic
set $\Gg^{inv,eq}\subset\Cc_G$ such that, for any $f\in\Gg^{inv,eq}$, if $x$
and $\tilde x$ are two equilibrium points of the ODE \eqref{eq} satisfying
$x_{i_0}=\tilde x_{i_0}$, then $x=\tilde x$. 
\end{prop}
\begin{demo}
First, notice that the case (c) is clearly implied by Corollary \ref{coro-eq}
since there is generically no equilibrium point in this case.

We use similar techniques as the ones in the proof of Proposition
\ref{prop-simpl}. Let $i_0$ be an observation cell and let $I$ be a set of
cells, which includes $i_0$ but not all the cells. We
denote by $\Gg^{inv,eq}_{I}$ the set of functions $f$ such that if $x$ and
$\tilde x$ are two equilibrium points of the ODE \eqref{eq} satisfying
$x_I=\tilde x_I$, then $x_i=\tilde x_i$ for at least one cell $i\not\in I$. 
We claim that $\Gg^{inv,eq}=\cap_{I} \Gg^{inv,eq}_{I}$ and thus that it is
sufficient to show that each $\Gg^{inv,eq}_{I}$ is a generic set. Indeed,
assume
that $f\in \cap_{I} \Gg^{inv,eq}_{I}$ and that $x$ and $\tilde
x$ are two equilibrium points such that $x_{i_0}=\tilde x_{i_0}$. Let $I$ be
the set of
cells $i$ such that $x_i=\tilde x_i$. Notice that $i_0\in I$ and assume that $I$
is not the whole set of cells. Since $f\in \Gg^{inv,eq}_{I}$, we
should have that $x_i=\tilde x_i$ for at least one cell $i\not\in I$, which is
obviously in contradiction with the definition of $I$. Thus, $I=\{1,\ldots,N\}$
and $x=\tilde x$.

We set $J$ to be the set of cells such that all their direct inputs belong to
$I$. We denote by $I^c$ and $J^c$ the complementary sets of cells. We use
Theorem \ref{th-Sard-Smale} with the following setting: 
$$\Mc=\{(x,\tilde x)\in X^2~,~x_I=\tilde x_I\text{ and }\forall
j\in I^c, x_j\neq \tilde x_j\}~,$$
$$\Nc=\RR^d\times\RR^{d_{J^c}}~,~~y=0\in\Nc~\text{ and }~\Lambda=\Cc_G~$$
We consider the function
$$\Phi:~(x,\tilde
x,f)\in\Mc\times\Lambda~\longmapsto~(\,f(x)\,,\,f_{J^c}(\tilde x)\,)\,
\in\RR^d\times\RR^{d_{J^c}}~.
$$
Let $(x,\tilde x,f)\in\Phi^{-1}(0)$.  We have $f(x)=0$ and $f_j(\tilde x)=0$ for
any $j\not\in J$. If $j\in J$, any direct input of $j$ belongs to $I$ and
therefore $f_j(\tilde x)=f_j(x)=0$. Hence, $x$ and $\tilde x$ are two
equilibrium points, which coincide on the cells $I$ but
differ in any cell in $I^c$. We claim that Hypothesis (ii) of Theorem
\ref{th-Sard-Smale} holds. Indeed, the map 
$$g\in\Cc_G \longmapsto D_f\Phi(x,\tilde x,f).g=(g(x),g_{J^c}(\tilde x))$$
is onto $\RR^{d+d_{J^c}}$ because, for any $j\in {J^c}$, $\tilde x$ differs
from $x$ in at least one direct input of $j$ and thus the value of $g_j(x)$ and
$g_j(\tilde x)$ can be chosen independently. Assume for a moment that
$dim(\Mc)<dim(\Nc)$, then Theorem \ref{th-Sard-Smale} implies that, generically
with respect to $f$, there are no equilibrium points $x$ and $\tilde x$ which
coincide on $I$ but differ on $I^c$, i.e. if $x$ and  $\tilde
x$ coincide on $I$, they must also coincide in at least one other cell. This
proves that $\Gg^{inv,eq}_{I}$ is a generic set.

It remains to show that, if the graph $G$ satisfies either (a) or (b), then we
have $dim(\Mc)<dim(\Nc)$, that is that the dimension $d_I$ is larger than the
dimension $d_J$. In the case of property (b), this trivially follows from the
definition of dimensionally decreasing graph since $I$ contains all the direct
input of $J$. Assume that $G$ is self-dependent. Then, any cell is a direct
input of itself and thus $J\subset I$. Therefore, $d_I=d_J$ if and only if
$J=I$. But, this would mean that any direct input of $J$ belongs to $J$ and
thus any indirect input of $I=J$ belongs to $I$. Since we assumed that $I$ is
not the whole set of cells, that $I$ contains $i_0$ and that $i_0$ is an
observation cell, this is not possible and $d_I>d_J$. 
\end{demo}

We believe that Proposition \ref{prop-eq-diff} is optimal that is
that if $G$ satisfies neither (a), (b) or (c), then one can find a circular
sub-network in which a counter-example similar to \eqref{ce-eq} can be
constructed.


\section{Generic inverse problem}\label{sect-inv}
The purpose of this section is to prove the inverse problem stated in Theorem
\ref{th1}. Its proof is split into two propositions. 

\begin{prop}\label{prop-inv-1}
Let $i$ be a cell, let $\tau$ be a positive number.
There exists a generic subset $\Gg^{traj,1}_{i,\tau}$ of $\Cc_G$ such
that the following property holds for any $f\in\Gg^{traj,1}_{i,\tau}$. Let
$x(t)$ and $\tilde x(t)$ be two solutions of \eqref{eq}. Assume that
$x_i(t)=\tilde x_i(t)$ on the time interval $[0,\tau]$ and assume that at least
one cell $j$ being a direct input of $i$ is such that $x_j(t)$ is not
constant on $[0,\tau]$. Then, for any cell $k$ which is a direct input of $i$,
the curves $\{x_k(t),t\in[0,\tau]\}$ and $\{\tilde x_k(t),t\in [0,\tau]\}$ must
cross.
\end{prop}
\begin{demo}
The proof of Proposition \ref{prop-inv-1} uses Theorem \ref{th-Henry} in the
same way as the proof of Proposition \ref{prop-eq}.  
The set of parameters is $\Lambda=\Cc_G$. 
We set $\Nc=(\Cc^0([0,\tau],\RR^d))^2$ and $y=(0,0)$. Let $j$ and $k$ be two
(possibly equal) direct input cells of the cell $i$. We set 
\begin{align*}
\Mc=\{&(\gamma,\tilde\gamma)\in\Cc^1([0,\tau],X)^2~,~\forall
t\in[0,\tau],~\gamma_i(t)=\tilde
\gamma_i(t),\\
&~\forall (s,t)\in[0,\tau]^2,~\gamma_k(t)\neq\tilde\gamma_k(s)\text{ and
}t\mapsto \gamma_j(t)\text{ is not constant }\}~,
\end{align*}
We define the function $\Phi\in\Cc^1(\Mc\times\Lambda,\Nc)$ by
$$\Phi(\gamma,\tilde\gamma,f)=\left(\begin{array}{l}
\frac{d}{dt}\gamma(t)-f(\gamma(t))\\
\frac{d}{dt}\tilde\gamma(t)-f(\tilde\gamma(t))\\
\end{array}\right)~.$$
The tangent space $T_\gamma\Mc$ is given by the
$(\omega,\tilde\omega) \in(\Cc^1([0,\tau],\RR^d))^2$ such that
$\omega_i\equiv\tilde\omega_i$. We have
$$D\Phi(\gamma,f).(\omega,g)=\left(\begin{array}{l}
\frac{d}{dt}
\omega(t)-Df(\gamma(t)).\omega(t)-g(\gamma(t))\\
\frac{d}{dt}
\tilde\omega(t)-Df(\tilde\gamma(t)).\tilde\omega(t)-
g(\tilde\gamma(t))
\end{array}\right)~.$$

{\noindent\it Step 1: $D_{(\gamma,\tilde\gamma)}\Phi$ is a left-Fredholm map.}\\
The arguments are the same as in the proof of Proposition \ref{prop-eq}. We use
a decomposition $D_{(\gamma,\tilde\gamma)}\Phi=L+K$ and apply Proposition
\ref{prop-Fred-1}. The map $L$ is given by $(\omega,\tilde \omega) \in
T_\gamma\Mc \mapsto (\frac{d}{dt} \omega,\frac{d}{dt}\tilde\omega) \in
(\Cc^0([0,\tau],\RR^d))^2$. Its kernel is the set of the constant functions such
that
$\omega_i=\tilde\omega_i$. Its image $R(L)$ consists in the functions
$(\theta,\tilde\theta)\in (\Cc^0([0,\tau],\RR^d))^2$ such that
$\theta_i\equiv\tilde\theta_i$. A closed complementary space $Y_2$ of $R(L)$ is
given by the space of the functions $(\theta,0)\in (\Cc^0([0,\tau],\RR^d))^2$
such that $\theta_k\equiv 0$ except if $k=i$.

{\noindent\it Step 2: checking Hypothesis (ii) of Theorem \ref{th-Henry}}\\
Once again, the arguments are similar to the ones of the proof of Proposition
\ref{prop-eq}. We have seen that $D_{(\gamma,\tilde\gamma)}\Phi$ can be written
$L+K$ as in Proposition \ref{prop-Fred-3}. It is sufficient to 
exhibit a sequence of functions $(g^n)\in \Cc_G$ such that $(D_f\Phi.g^n)$
composes a free family of $Y_2$. By assumption, the curve
$\Gamma=\{(\gamma_j(t),\gamma_k(t)),t\in[0,
\tau ] \} \subset X_j\times X_k$ is not reduced to a point and does not meet
the curve
$\tilde\Gamma=\{(\tilde\gamma_j(t),\tilde\gamma_k(t)), t\in [ 0 , \tau]\}$ (if
$j=k$ simply consider $\Gamma=\{\gamma_j(t)\}$ and
$\tilde\Gamma=\{\tilde \gamma_j(t)\}$).
We construct the suitable vector fields $g^n$ exactly as in the proof of
Proposition \ref{prop-eq}, introducing disjoints open sets $\Oc^n_{j,k}\subset
X_j\times X_k$ which intersect the curve $\Gamma$ but not the curve $\tilde
\Gamma$.

{\noindent\it Step 3: $\Mc$ and $\Lambda$ are obviously separable}

{\noindent\it Step 4: conclusion}\\
We can apply Theorem \ref{th-Henry}, which shows that there exists a generic
subset of admissible vector fields $f\in\Cc_G$ such that $0$ is not in the
image of $\Phi(.,f)$. This means that there are no solutions
$\gamma(t)$ and $\tilde \gamma(t)$ of \eqref{eq} on the time interval
$[0,\tau]$, being
equal on the cell $i$, and such that $x_j(t)$ is not constant on $[0,\tau]$ and
the curves $\{\gamma_k(t),t\in[0,\tau]\}$ and $\{\tilde \gamma_k(t),t\in
[0,\tau]\}$ do not intersect.
\end{demo}

In Proposition \ref{prop-inv-1}, we had to exclude the case of constant direct
inputs, because we wanted to construct an infinite-dimensional complementary
subspace in the second step of the proof. The case of constant direct inputs
follows from Section \ref{sect-eq}.
\begin{prop}\label{prop-inv-2}
Let $i$ be a given cell, let $I$ be the set of all its indirect inputs and
assume that $i$ does not belongs to $I$. Then, there exists a dense open set
$\Gg^{traj,2}_i$ of $\Cc_G$ such that, for any $f\in \Gg^{traj,2}_{i}$, the
following property holds. Let $x(t)$ and $\tilde x(t)$ be two trajectories of
\eqref{eq} such that $x_I(t)$ and $\tilde x_I(t)$ are constant. If
$x_i(0)=\tilde x_i(0)$ and $f_i(x(0))=f_i(\tilde x(0))$ then $x_j(0)=\tilde
x_j(0)$ for any $j$ being a direct input of $i$.
\end{prop}
\begin{demo}
The set of cells $I$ is an independent sub-network that is that 
the ODE \eqref{eq} restricted to this set is another coupled cell network 
\begin{equation}\label{eqI}
\dot x_I(t)=f_I(x_I(t))~.
\end{equation}
Applying Proposition \ref{prop-simpl} to this
network, we get the existence of a dense open set $\Gg^{simpl}_I\subset\Cc_G$
such that for any $f\in \Gg^{simpl}_I$, there is only a finite number of
equilibrium points for \eqref{eqI}. Moreover, these equilibrium points are
simple and so depend smoothly on $f_I$. 

Now, since $i\not\in I$, the fact that $f\in
\Gg^{simpl}_I$ does not depend on $f_i$.
Therefore, we can easily
perturb $f_i$ so that $f_i(x_I)\neq f_i(\tilde x_I)$ for any couple of 
equilibrium points $(x_I,\tilde x_I)$ of \eqref{eqI}, which do not coincide
on the direct inputs of $i$ (otherwise of course $f_i(x_I)=f_i(\tilde
x_I)$). Moreover, this property holds also for small perturbations of $f$
since $f_i$ and the equilibria $x_I$ depend smoothly on $f$. Thus, we get a
dense open set $\Gg^{traj,2}_i$ of functions $f$ such that if $x$ and
$\tilde x$ are two trajectories, which are constant on $I$, then
$f_i(x)=f_i(\tilde x)$ if and only if $x_j=\tilde x_j$ for any $j$ being a
direct input of $i$.
\end{demo}

{\noindent\bf Proof of Theorem \ref{th1}: } let $\Gg^{traj}$ be the
intersection of all the previous generic subsets $\Gg^{simpl}$, $\Gg^{eq}$,
$\Gg^{inv,eq}$, $\Gg^{traj,1}_{i,1/n}$ and $\Gg^{traj,2}_i$, for all
cell $i$ and all $n\in\NN^*$. Notice that $\Gg^{traj}$ is a
generic subset of $\Cc_G$, since it is a countable intersection of generic
subsets. 

Assume that $f\in\Gg^{traj}$ and that $i_0$ is an observation cell. Let
$x(t)$ and $y(t)$ be two trajectories of \eqref{eq} such that
$x_{i_0}(t)=y_{i_0}(t)$ on some time interval $(a,b)$. Assume that $x$ and $y$
do not coincide in all the direct input cells of $i_0$, we would like to obtain
a contradiction. First, notice that, up to translating and restricting the
interval of times, we can replace $(a,b)$ by $[0,1/n]$ and assume that the
curves $\{x_k(t),t\in [0,1/n]\}$ and $\{y_k(t),t\in [0,1/n]\}$
do not cross for at least one of the direct inputs of $i_0$.\\
{\it First case: $x_j(t)$ or $y_j(t)$ is not constant in $[0,1/n]$ for at least
one of the direct inputs $j$ of $i_0$.} Since $f\in
\Gg^{traj,1}_{i_0,n}$ applying Proposition
\ref{prop-inv-1}, we obtain a contradiction.\\
{\it Second case: $x_j(t)$ and $y_j(t)$ are constant in $[0,1/n]$ for any
direct input $j$ of $i_0$, and $i_0$ is not an indirect input of itself.}
Since $f\in\Gg^{eq}$, Property (a) of Theorem \ref{th2} shows that $x(t)$ and
$y(t)$ are constant for any indirect input of $i_0$. Since
$f\in\Gg^{traj,2}_{i_0}$, Proposition \ref{prop-inv-2} shows that $x(t)$ and
$y(t)$ coincide in any direct input of $i_0$, which is in contradiction with our
assumption.\\
{\it Third case: $x_j(t)$ and $y_j(t)$ are constant in $[0,1/n]$ for any
direct input $j$ of $i_0$, and $i_0$ is an indirect input of itself.} Since
$f\in\Gg^{eq}$, Property (a) of Theorem \ref{th2} shows that $x(t)$ and
$y(t)$ are constant for any cell being in the set $I$ of the indirect inputs of
$i_0$. Notice that $i_0\in I$ and thus $x_I(t)$ and $y_I(t)$ are two
equilibrium points of the independent sub-network $I$, which coincide on the
observation cell $i_0\in I$. Applying Proposition \ref{prop-eq-diff}, we obtain
that $x(t)$ and $y(t)$ coincide in $I$ except if the graph $G$ restricted to
the set of cells $I$ is dimensionally non-increasing and is neither a
self-dependent directed graph nor a dimensionally decreasing graph. Notice that
the fact that $G$ restricted to $I$ is dimensionally non-increasing can be
omitted: if it is not the case, then Corollary \ref{coro-eq} shows that there is
no constant solution to \eqref{eqI} and thus $x_I(t)$ and
$y_I(t)$ cannot be two equilibrium points of the independent sub-network $I$.

By applying the above arguments recursively, we get that $x(t)$ and $y(t)$
coincide in any cell in some interval of times $[0,1/n]$, {\it except maybe if}
in our recursion process, we meet the exception of the above third case. This
exception may happen only if the trajectories $x(t)$ and
$y(t)$ are both constant in a set of cells $I$ and the graph $G$ restricted to
$I$ is neither a
self-dependent directed graph nor a dimensionally decreasing graph. To finish
the proof of Theorem \ref{th1}, notice that if $x(t)$ and $y(t)$
coincide in any cell in some time interval $[0,1/n]$, then they must
coincide  for any time $t$ since \eqref{eq} is a classical ODE.
\hfill $\square$


\section{Generic observability of oscillations and
stabilisation}\label{sect-coro}

In this section, we show that Properties (b) and (c) of Theorem \ref{th2} and
Theorem \ref{th3} are straightforward consequences of the inverse problems
studied in the previous sections.

\vspace{3mm}

{\noindent\bf Proof of Properties (b) and (c) of Theorem \ref{th2}: }
Let $\Gg^{osc}$ be the generic set $\Gg^{traj}$ introduced in the proof of
Theorem \ref{th1} and assume that $f\in\Gg^{osc}$. Let $x(t)$ be a solution of
\eqref{eq} and assume that there is a cell $i$ such that $x_i(t)$ is
$T-$periodic on $(a,a+T+\tau)$. Let $I$ be the set of indirect input cells of
$i$. We set $y(t)=x(t+T)$. We
have $x_i(t)=y_i(t)$ on $(a,a+\tau)$. Since $f\in\Gg^{traj}$, using Theorem
\ref{th1}, we obtain that $x_I(t)=y_I(t)$ except maybe
in the cells where $x$ and $y$ are constant. But, since $y(t)=x(t+T)$, if both
are constant in time, they must be equal. To conclude, we get that
$x_I(t)=y_I(t)=x_I(t+T)$ for all $t\in\RR$, and thus any cell in $I$ is
$T-$periodic, which proves Property (b).

Let us show Property (c). Notice that, in Property (b), it is possible
that an indirect input cell $j$ of $i$ is $T'-$periodic with $T'<T$. Therefore,
if $x_i(t)$ is exactly $T-$periodic, it is possible that $x_j(t)$ is not exactly
$T-$periodic but exactly $T'-$periodic with $T'<T$ (or even constant). However,
if $G$ is a strongly connected graph, then applying Property (b) (or (a)) to the
cell $j$, for which $i$ is an indirect input cell, we get that $x_i(t)$ is also 
exactly $T'-$periodic (or constant), which is of course impossible since $T\neq
T'$. Therefore, if $G$ is a strongly connected graph, then $x_j(t)$ must be
also exactly $T-$periodic.
\hfill $\square$

\vspace{3mm}

{\noindent\bf Proof of Theorem \ref{th3}: } Let $\Gg^{stab}$ be the generic
subset $\Gg^{traj}$ of $\Cc_G$ introduced in the proof of
Theorem \ref{th1} and assume that $f\in\Gg^{stab}$. Let $x(t)$ be a
solution of \eqref{eq} and let $i$ be an observation cell. Assume that $x_i(t)$
converges to a constant $x^*_i$ when $t$ goes to $+\infty$. The main ingredient
of this proof is the classical concept of $\omega-$limit set. We recall that
the $\omega-$limit set of $x(0)$ is defined as
\begin{align}
\omega(x(0))&=\cap_{n\in\NN} \overline{\{x(t)~,~t\in[n,+\infty)\}}\nonumber\\
&=\{ \xi\in X~,~\text{there exists a sequence }(t_n)\text{ converging to
}+\infty\text{ such that }\nonumber \\
&~~~~~~x(t_n)\text{ converges to }\xi\text{ when }n\text{ goes to }
+\infty\}~. \label{def-omega}
\end{align}
Since $X$ is a compact manifold, for all $x(0)\in X$, $\omega(x(0))$ is a
connected and compact non-empty set. Moreover $\omega(x(0))$ is invariant by
the dynamical system generated by the ODE \eqref{eq}, i.e. the $\omega-$limit
set is a union of complete trajectories of \eqref{eq} (see for example
\cite{Palis-de-Melo}). By continuity, we know that, in our particular case, any
trajectory $\xi(t)$ of
$\omega(x(0))$ satisfies $\xi_i(t)=x^*_i$. By using Property (a) of Theorem
\ref{th2}, we get that $\xi(t)$ is an equilibrium point of \eqref{eq}. Since
$f\in\Gg^{simpl}$, any equilibrium is simple and hence isolated. Since
$\omega(x(0))$ is connected, there is exactly one equilibrium point $\xi=x^*$
in $\omega(x(0))$.  Since $X$ is compact, this shows that
$x(t)$ converges to $x^*$.
\hfill $\square$


\section{Discussion on further results}\label{sect-further}

In this section, we state several adaptations or generalisations of our main
results. They are presented as claims and their proofs are only
roughly sketched. By ``claims'', we mean that the
results hold, but no complete proofs are given. The proofs are
mutatis mutandis the same as the ones of Theorem \ref{th1}, \ref{th2} or
\ref{th3}, we only give short arguments, that are, we hope, sufficiently
convincing.

\vspace{5mm}

{\noindent \bf The non-compact case.}\\
It is possible to adapt the main results of this paper to the non-compact case.
Here, we typically choose $X_i=\RR^{d_i}$ to fix the notations. The first
problem is to define
a suitable topology on the space $\Cc_G$ of admissible vector fields, i.e. a
topology on $\Xg^1(X)$. The classical strong $\Cc^1-$topology is too
restrictive since it applies only on bounded functions. Therefore, one usually
endow $\Xg^1(X)$ with the Whitney topology, that is the topology generated by
the neighbourhoods
$$
\Oc_{f,\delta}=\{g\in\Xg^1(X)~/~\max(\|f(x)-g(x)\|,\|Df(x)-Dg(x)\|)<\delta(x),
~\forall x\in X\}~,
$$
where $f$ is any function in $\Xg^1(X)$ and where $\delta$ is any continuous
positive real function on $X$. We refer to \cite{GG} for a
complete discussion on this topology. The main issue is that $\Xg^1(X)$
endowed with the Whitney topology is not a metrizable space and the close sets
are not the sequentially closed sets. In particular, we cannot directly use
Theorems \ref{th-Sard-Smale} or \ref{th-Henry} with $\Lambda=\Cc_G$. However,
it is noteworthy that $\Cc_G$ endowed with the Whitney topology is still a
Baire space and thus generic subsets are dense subsets.

A second issue is the global existence of solutions of the differential
equation \eqref{eq}. In order to ensure it, we restrict the class of vector
fields to
$$\tilde\Cc_G=\{~f\in\Cc_G~,~\exists M>0,~\forall x\in X,~\|x\|\geq
M\Rightarrow \langle f(x)|x\rangle <0~\}~.$$
If $f\in\tilde\Cc_G$, any solution $x(t)$ of \eqref{eq} exists for all $t\geq
0$ and eventually belongs to the ball $B_X(0,M)$. In particular, the
$\omega-$limit set defined by \eqref{def-omega} is
non-empty.
\begin{claim}\label{claim1}
There exists a generic subset $\tilde\Gg$ of $\tilde \Cc_G$ such that, for any
$f\in\tilde\Gg$, all the observation properties stated in Theorems \ref{th1},
\ref{th2} and \ref{th3} hold. 
\end{claim}
\begin{sketch}
Since we cannot apply Theorems \ref{th-Sard-Smale} or \ref{th-Henry} with
$\Lambda=\tilde\Cc_G$, we use the following strategy, which has now become
classical:\\
- split the manifold $X$ in compact subsets,\\
- consider the desired property in each of these compact subsets and show
that it is open,\\
- prove the density by applying a transversality theorem in an
open subset of $X$ containing the compact subset.\\
To give an example, let us explain how we adapt Proposition
\ref{prop-simpl}. We cover $X$ by a countable union of closed balls
$\overline B_X(0,n)$. We introduce the set $\tilde \Gg_n^{simpl}$ consisting in
functions $f\in\tilde\Cc_G$ such that any equilibrium point $e$ of \eqref{eq}
satisfying $e\in \overline B_X(0,n)$ is simple. Arguing as in the proof of
Proposition \ref{prop-simpl}, we show that $\tilde \Gg_n^{simpl}$ is open in
$\tilde\Cc_G$. To be able to use Theorem \ref{th-Sard-Smale}, we introduce the
space $\Lambda$ of admissible vector fields $g\in\Cc_G$ being compactly
supported in $\overline B_X(0,n+1)$. Notice that $\Lambda$ is a Banach space
endowed with the classical strong $\Cc^1-$topology. Then, we set
$\Mc=B_X(0,n+1)$, $\Nc=\RR^d$, $y=0\in\RR^d$ and 
$$\Phi:~(x,g)\in\Mc\times\Lambda~\longmapsto~(f+g)(x)\in\RR^d~.$$
Using the same arguments as in the proof of Proposition
\ref{prop-simpl}, we show that there exists a generic set of $\Cc_G$ such that
any equilibrium $e$ of \eqref{eq} with $\|e\|<n+1$ is simple. Then notice that,
if $g$ is close to zero in $\Lambda$, then $f+g$ is close to $f$ in
$\tilde\Cc_G$. As a consequence, we can perturb $f$ such that any
equilibrium $e$ of \eqref{eq} with $\|e\|<n+1$ is simple, and a fortiori
this shows that $\tilde \Gg_n^{simpl}$ is dense. Finally
$\tilde\Gg^{simpl}=\cap_n \tilde \Gg_n^{simpl}$ is a generic subset of
$\tilde\Cc_G$.
\end{sketch}

\vspace{5mm}

{\noindent \bf General manifolds.}\\
Assume that for each cell $i$, $X_i$ is a compact $\Cc^2-$manifold of dimension
$d_i$. 

The first problem occurs when one would like to define a network with a
graph which is not self-dependent. Indeed, it may be abusive to assume that a
vector field $f_i(x)$ in a cell $i$ does not depend on $x_i$ since $f_i(x)\in
T_{x_i}X_i$. In fact, to define a vector field $f_i$ which is
independent of $i$, one has to be able to introduce a notion of constant vector
fields on $X_i$. In the cases of $X_i=\TT^{d_i}$ or of $X_i=\RR^{d_i}$ as above,
we have identified each $T_{x_i}X_i$ to $\RR^{d_i}$, which enables us to
define a notion of constant vector field (even if we have not explicitely
mentioned this, the concept is natural). If $X_i$ is a
$2n$-dimensional sphere, then any concept of constant vector field will be
irrelevant since any continuous vector field must vanish. A suitable notion of
constant vector fields would be to assume that $X_i$ is a parallelizable
manifold, i.e. that there exist $d_i$ vector fields $\xi_1$, $\ldots$,
$\xi_{d_i}$ such that for each $x_i\in X_i$, $(\xi_k(x_i))$ is a basis of
$T_{x_i}X_i$, and then to define constant vector fields as linear combinations
of the $\xi_k$. The spaces $\RR^{d_i}$, the tori, the sphere
$S^3$ and their products are parallelizable manifolds.
However, except these ones, there is no other natural examples. For this
reason, if one deals with general manifolds, it is natural to assume that
the graph $G$ is self-dependent.

Except for the problem of self-dependency, our results can be extended to
general compact manifolds. 
\begin{claim}\label{claim2}
Assume that each $X_i$ is a compact $\Cc^2-$manifold and assume that $G$ is a
self-dependent graph. Then, Theorems \ref{th1}, \ref{th2} and \ref{th3} hold.
\end{claim}
\begin{sketch}
When trying to repeat the proofs of these theorems in the framework of general
manifolds, the main difficulty is to deal with global abstract setting for
manifolds. The following questions arise. How to define abstractly and properly
the Banach manifold of
$\Cc^1-$ paths on $X$? How to check the surjectivity of a functional $D\Phi$,
whose image is now included in the space $T(TX)$? 
A simple way to overcome these difficulties is to work in local coordinates.
We consider the canonical product manifold structure for
$X=X_1\times\ldots\times X_N$. In this way the charts are not mixing the
coordinates of the different cells of the network. Let us simply give an example
of how we have to modify our propositions. The statement of Proposition
\ref{prop-eq} becomes
\begin{quote}
Let $i$ and $j$ be two cells such that $j$ is a direct input of $i$. Let
$\tau>0$ be given and let $(U,\alpha)$ be a chart of an atlas $\Uc$. Let
$\Gg^{eq}_{i,j,\tau,U}$ be the set of all the
admissible vector fields $f\in\Cc_G$ such that, if $x(t)$ is a solution of
\eqref{eq} such that $x_i(t)$ belongs to $U$ for all $t\in[0,\tau]$ and is
constant on this interval, then $x_j(t)$ is also constant on
$[0,\tau]$. Then, $\Gg^{eq}_{i,j,\tau,U}$ is an generic subset of $\Cc_G$.
\end{quote}
The proof is the same as the one of Proposition \ref{prop-eq}, 
except that we now work in the local chart $U$ and consider the space
$$\Mc=\{\gamma \in\Cc^1([0,\tau],\alpha(U)),~\gamma_i(t)\text{ is constant and
}\gamma_j(t)\text{ is not constant }\}$$
and the function
$$\Phi(\gamma,f)=\frac{d}{dt}\gamma(t)-D\alpha(\alpha^{-1}(\gamma(t))).f(\alpha^
{ -1}(\gamma(t)))~.$$
Since the manifolds $X_i$ are covered by a finite number of charts, we recover
our global results by intersecting the sets $\Gg^{eq}_{i,j,\tau,U}$, which have
been locally obtained.
\end{sketch}

Of course, we could even state our results in the framework of
non-compact general manifolds by mixing Claims \ref{claim1} and \ref{claim2}.

\vspace{5mm}

{\noindent \bf Prevalent results.}\\
In this paper, we used the genericity to give a meaning to the notion
of ``almost everywhere''. It is noteworthy that it is not the only notion of
``almost everywhere'' in Banach spaces. Except the habits, there is no
real reason to give more importance to genericity in comparison to other
natural notions of ``almost everywhere''. In particular, the prevalence has
recently attracted much attention.

Let $X$ be a Banach space. Christensen in \cite{Christ} introduced the notion
of Haar-nul set: a Borel set $B$ of $X$ is said \emph{Haar-nul} if
there exists a finite non-negative measure $\mu\not\equiv 0$ with compact
support such that for all $x\in X$, $\mu(x+B)=0$. More generally, any set
$B\subset X$ is said Haar-nul if it is contained in a Haar-nul Borel set.
Let $U$ be an open subset of $X$. A set $P\subset U$ is said \emph{prevalent} in
$U$ if $U\setminus P$ is a Haar-nul set of $X$. The notion of prevalence first
appeared in the work of Hunt, Sauer and Yorke \cite{HSY}, where it is proved
that prevalent sets are dense and that a countable intersection of prevalent
sets is prevalent. It is interesting to see that, in finite-dimensional spaces,
prevalent sets are exactly the sets of full Lebesgue measure (see \cite{HSY}).
We refer to \cite{Ott-Yorke} for a review on prevalence. 

We proved in \cite{RJ-prevalence} that, if $\Lambda$ is an open subset of a 
Banach space, then in the transversality theorems, Theorems \ref{th-Sard-Smale}
and \ref{th-Henry}, the set $\Gg$ is not only generic but also prevalent. Thus,
we obtain the following generalisation of our main results.
\begin{claim}
Any generic subset appearing in Theorems \ref{th1}, \ref{th2} and \ref{th3}, as
well as in the other results of this paper, is not only a generic subset of
$\Cc_G$ but also a prevalent subset. 
\end{claim}

\vspace{5mm}

{\noindent \bf The discrete-time case: cellular automata.}\\
Instead of considering the time-continuous dynamical system generated by the ODE
\eqref{eq}, we could also study the discrete-time model
\begin{equation}\label{eq-disc}
x^0\in X~,~~x^{n+1}=f(x^n)~~~~~f\in\Cc_G~.
\end{equation}
The discrete-time model is also meaningful for applications. If the cells form
a grid and if the direct inputs are exactly the neighbouring cells of each
cell, then we get a usual cellular automaton. The observation problems stated
in this paper are obviously transposed into the discrete-time model.
\begin{claim}
Let $G$ be a given graph. There exists a time $K\in\NN$ depending on $G$ and
the dimension $d$ of the phase space $X$ such that:\\
- Theorem \ref{th3} holds,\\
- Theorem \ref{th1} holds as soon as the time interval $(a,b)$ satisfies
$b-a\geq K$,\\
- Theorem \ref{th2} holds as soon as $b-a\geq K$ and $\tau\geq K$.
\end{claim}
\begin{sketch}
In the continuous-time model, as soon as we have considered a non-constant
trajectory in some interval of times, we get an infinite-dimensional space of
freedom, as one can see in the proof of Proposition \ref{prop-eq} for
example. There, we may have used Theorem \ref{th-Henry} with $Z$
infinite-dimensional in Hypothesis (ii). The difficulty with the discrete-time
model is that this is no longer possible. Exactly as in the proof of Proposition
\ref{prop-eq-diff}, we have to be careful in counting the dimension. By using
techniques similar to the ones of Proposition \ref{prop-eq-diff}, we can
however adapt the results of this paper into the frame of the discrete-time
model \eqref{eq-disc}. The time $K$ is basically equal to $dD$ where $d$ is the
dimension of $X$ and $D$ is the diameter of the graph, i.e. the larger distance
between one cell and one of its indirect inputs.
\end{sketch}

\vspace{5mm}

{\noindent \bf Generic hyperbolicity of equilibrium points.}\\
In Proposition \ref{prop-simpl}, we proved the generic simplicity of
equilibrium points. We could wish to go further by proving the generic
hyperbolicity of equilibrium points. We recall that an equilibrium $e$ is
hyperbolic if $Df(e)$ has no spectrum on the vertical line $i\RR$. The
hyperbolicity of an equilibrium point $e$ implies that the flow near $e$ is
qualitatively the same as the one of the linear equation $\dot
x(t)=Df(e).x(t)$, see \cite{Palis-de-Melo}. In classical ODEs, the hyperbolicity
of equilibrium points is generic. However, it is not always the case in coupled
cell networks. Indeed, consider the simple counter-example \eqref{ce-eq} with
$X_1=X_2=\RR/\ZZ$. It is not difficult to construct $f$ such that
there exists an equilibrium $e=(e_1,e_2)$, which is simple and satisfies
$f'_1(e_2)f'_2(e_1)<0$. Notice that the existence of such an equilibrium holds
in a small neighbourhood of $f$. Since 
$$Df(e)=\left(\begin{array}{cc} 0 & f'_1(e_2)\\f'_2(e_1) &0 
\end{array}\right)~$$
has $\pm i\sqrt{|f'_1(e_2)f'_2(e_1)|}$ for eigenvalues, 
such an equilibrium is never hyperbolic.

Fortunately, we can at least show the generic hyperbolicity of equilibrium
points in the meaningful case of self-dependent graph.
\begin{claim}\label{cl-1}
Assume that $G$ is a self-dependent graph. Then, there exists an open dense set
$\Gg^{hyp}$ of $\Cc_G$ such that for any $f\in\Gg^{hyp}$, all the equilibrium
points of \eqref{eq} are hyperbolic.
\end{claim}
\begin{sketch}
We recall that a simple
(resp. hyperbolic) equilibrium depends smoothly of $f$ and remains simple (resp.
hyperbolic) for small perturbations of $f$. Thus, the set $\Gg^{hyp}$ is
open and it remains to prove its density. Let $f\in\Cc_G$. Up to perturb it, we
can assume that $f$ belongs to $\Gg^{simpl}$, the open dense set of Proposition
\ref{prop-simpl}. Then, there is at most a finite number of equilibrium points
of \eqref{eq}. It is sufficient to
understand how to make a simple equilibrium point $e$ hyperbolic by perturbing
$f$, since then we can apply several times the same type of perturbations to
make, one by one, all the equilibria of $f$ hyperbolic. Notice that, if $e$
is a simple equilibrium, then $e$ is an hyperbolic equilibrium of \eqref{eq} for
$\tilde f$ being close to $f$ and satisfying $\tilde f(e)=f(e)=0$ and $D\tilde
f=Df+\varepsilon Id$. The final central argument is that such a perturbation of
$f$, consisting in perturbing the differential by a homothety, is
possible in the class $\Cc_G$ of admissible vector fields if and
only if $G$ is self-dependent.
\end{sketch}


\end{document}